\documentclass[11pt,a4paper,reqno]{amsart}
\usepackage{amsmath}
\usepackage{amsthm}
\usepackage{amssymb}
\usepackage{graphics}
\usepackage{graphicx}
\usepackage{colortab} 
\usepackage{enumitem}

\usepackage{colortab}

\usepackage[colorlinks=true,linkcolor=blue,citecolor=blue]{hyperref}
\usepackage{prettyref}

\theoremstyle{plain}
  \newtheorem{theorem}{\bf Theorem}[section]
  \newtheorem{proposition}[theorem]{\bf Proposition}
  \newtheorem{lemma}[theorem]{\bf Lemma}
  \newtheorem{corolario}[theorem]{\bf Corollary}

\theoremstyle{remark}
  \newtheorem{remark}[theorem]{\bf Remark}
  \newtheorem{definition}[theorem]{\bf Definition}
 
	 \newcommand{\R}{\mathbb{R}}

	 \newcommand{\N}{\mathbb{N}}
	 \newcommand{\calK}{\mathcal{K}}

	 \def\sech{sech}
	 \def\distr{C^{\infty} (\mathbb{S}^{d-1})^*}

\begin{document}
\title[The Fourier extension operator of distributions in $H^\alpha(\mathbb{S}^{d-1})$]{The Fourier extension operator of distributions in Sobolev spaces of the sphere and the Helmholtz equation}
\author{ J. A. Barcel\'{o}, M. Folch-Gabayet, T. Luque, S. P\'{e}rez-Esteva, 
    M. C. Vilela}\thanks{
       The first and the fifth authors are supported by the Spanish Grant MTM2017-85934-C3-3-P. The second and the fourth authors are supported by the Mexican proyect PAPIIT-UNAM IN106418. The third author is supported by the spanish grant MTM2017-82160-C2-1-P.}


\maketitle

\begin{abstract}
The purpose of this paper is to characterize all the entire solutions of the homogeneous Helmholtz equation (solutions in $\R^d$) arising from the Fourier extension operator of distributions in  Sobolev spaces of the sphere  $H^\alpha(\mathbb{S}^{d-1}),$ with $\alpha\in \mathbb{R}$. 
We present two characterizations. 
The first one is written in terms of certain $L^2$-weighted norms involving real powers of the spherical Laplacian. 
The second one  is in the spirit of the classical description of the Herglotz wave functions given by P. Hartman and C. Wilcox. For $\alpha>0$ this characterization involves  a multivariable square function evaluated in a vector of entire solutions of the Helmholtz equation, while for  $\alpha<0$ it is written in terms of an spherical integral operator acting as a fractional integration operator.
Finally, we also characterize all the solutions that are the Fourier extension operator of  distributions in the sphere. 
\end{abstract}

	\section{Introduction and statement of results}
Consider the  Fourier extension operator
\begin{equation}\label{def:extensionOp}
\mathcal{E}\phi(x):=\widehat{\phi d\sigma}(x)=\int_{\mathbb{S}^{d-1}}e^{-ix\cdot\xi}\phi(\xi)\,d\sigma(\xi),\quad \phi\in L^1(\mathbb{S}^{d-1}),\ x\in\mathbb{R}^d,
\end{equation}

where $\mathbb{S}^{d-1}$ is the unit sphere in $\mathbb{R}^d$, $d\sigma$ is the surface measure on $\mathbb{S}^{d-1}$ and $\phi d\sigma $ is the temperate distribution (a singular Borel measure in $\R^d$) defined by $\sigma$ with density $\phi$. This very important operator in harmonic analysis is a source of entire solutions  of the homogeneous Helmholtz equation; namely, $u=\mathcal{E}\phi$ satisfies the  equation
$$\Delta u+u=0,$$
 in $\mathbb{R}^d$.
Specially important is the case  when  $\phi \in L^2(\mathbb{S}^{d-1})$.
 In this case  the function $u=\mathcal{E}\phi$ is called  Herglotz wave function with density $\phi$. The space of Herglotz wave functions, that will be denoted by $W(\mathbb{R}^d)$, plays an important role in in the study of wave scattering problems. For more details on this, see the monograph \cite[pp.56-49 and pp.223-230]{CK}.

Several characterizations of $W(\mathbb{R}^d)$ are known. The first one was given by Hartman and Wilcox \cite{HartmanWilcox}, who proved that the Herglotz wave functions are precisely all the entire solutions $u$ of the homogeneous Helmholtz equation satisfying  the condition 
\begin{equation*}
\limsup_{R\rightarrow\infty}\,\frac{1}{R}
\int_{|x|<R}|u(x)|^2\,dx
<\infty.
\end{equation*} 
Furthermore, they proved that such condition can be replaced by
\begin{equation}\label{limite_herglotz}
\|u\|_{L}^{2}:=\lim_{R \rightarrow \infty} \frac{1}{R} \int_{|x|<R}|u(x)|^2 dx<\infty,
\end{equation}
and  hence by 
\begin{equation}\label{HartmanWilcox}
\|u\|_{A}^{2}:=\sup_{R>0}\,\frac{1}{R}
\int_{|x|<R}|u(x)|^2\,dx
<\infty.
\end{equation} 
The norm given by \eqref{limite_herglotz} is indeed a Hilbert space norm on $W(\mathbb{R}^d)$. Moreover, a Herglotz wave function $u$ with density $\phi$ satisfies that (see \cite[Theorem 2.2]{HartmanWilcox} and \cite[Theorem 3.30, for d=3]{CK})
\begin{align}\label{norm_equivalente}
\|u\|_{L}\sim\|u\|_{A}\sim\|\phi\|_{L^2(\mathbb{S}^{d-1})}.
\end{align}

More recently, a different characterization was given in terms of a weighted Sobolev norm. More precisely,  Herglotz wave functions  are all the entire solutions $u$ of the homogeneous Helmholtz equation such that 
\begin{equation}
\label{salvador}
\|u\|^2_S:=\int_{\mathbb{R}^d}\left(|u(x)|^2+\left|\nabla_S u(x)\right|^2\right)\langle x\rangle^{-3}\,dx
<\infty,
\end{equation}
where  $ \nabla_S u$ denotes the spherical gradient of $u$ defined by 
\begin{equation*}\label{spherical gradient}
\nabla_S u(x)= r \left(  \nabla u(x)- \frac{x}{r}\,\frac{\partial u}{\partial r}(x)  \right),\qquad r=|x|,
\end{equation*}
and $\langle x\rangle:=(1+|x|^2)^{1/2}$. 
The proof of this result for $d=2$ can be found in \cite[Theorem 1]{A+F-G+P-E}, and for $d>2$ in \cite[Theorem 4]{P-E+V-D} (see also \cite{BBR}). This new description of $W(\mathbb{R}^d)$ gives a clearer look of this space as a Hilbert space and its structure of a Hilbert space with reproducing kernel. Moreover, they prove that a Herglotz wave function $u$ with density $\phi$ satisfies that
\begin{align}\label{norm_equivalente2}
\|u\|_S\sim\|\phi\|_{L^2(\mathbb{S}^{d-1})}.
\end{align}

The extension operator given in \eqref{def:extensionOp} can be naturally defined in  $C^{\infty}(\mathbb{S}^{d-1})^*$, the space of  distributions on the sphere. Indeed, if $\mathcal{S}'(\R^{d})$ denotes the space of tempered distributions on $\R^d$, from the inclusion
\begin{equation*}
 C^{\infty}(\mathbb{S}^{d-1})^*\hookrightarrow\mathcal{S}'(\R^{d})
 \end{equation*} 
 given by 
 \begin{equation*}
 f(\varphi)=\langle f,\varphi_{\vert\mathbb{S}^{d-1}}\rangle, \quad\varphi\in\mathcal{S}(\R^d),
 \end{equation*}
  where $\langle\cdot,\cdot\rangle$ denotes the duality in $C^{\infty}(\mathbb{S}^{d-1})$,
 we can define for $f\in C^{\infty}(\mathbb{S}^{d-1})^*$, 
\begin{equation}
\label{definicionExtension}
 \mathcal{E}f=\mathcal{F}(f),
 \end{equation} 
where $\mathcal{F}$ denotes the Fourier transform. Moreover, since every $f\in C^{\infty}(\mathbb{S}^{d-1})^* $ defines a distribution with  compact support contained in $\mathbb{S}^{d-1}$, then $\mathcal{E}f$ is a $C^\infty $ function in $\R^d$ given by
 \begin{equation*}
   \mathcal{E}f(x)=\langle f,e^{-ix\cdot(\cdot)}\rangle,\quad x\in\mathbb{R}^d,
 \end{equation*}  
and thus,  $\mathcal{E}f$ is an entire solution of the homogeneous Helmholtz equation. 
 
The aim of this paper is to characterize all the entire solutions of the homogeneous Helmholtz equation arising from the extension operator of distributions in Sobolev spaces of the sphere  $H^\alpha(\mathbb{S}^{d-1}),$ with $\alpha\in \mathbb{R}$ (for a precise definition of  these Sobolev spaces, see below the Subsection \ref{preliminares}). The aimed characterizations will be in the sprit of  those given in \eqref{HartmanWilcox} and in \eqref{salvador} for the  particular case of Herglotz wave functions.  

\begin{definition}
	Let $\alpha\in\mathbb{R}$. We say that the function $u:\mathbb{R}^d\rightarrow\mathbb{C}$ is a $\alpha$-Herglotz wave function if $u=\mathcal{E}f$ for some $f\in H^\alpha(\mathbb{S}^{d-1})$, and we will denote by $\mathcal{W}^\alpha(\mathbb{R}^d)$ the space of all $\alpha$-Herglotz wave functions.
\end{definition}	
Notice that we recover the space of Herglotz wave functions when $\alpha=0$, that is, 
$\mathcal{W}^0(\mathbb{R}^{d})=W(\mathbb{R}^d)$. 

Characterizations \eqref{norm_equivalente} and \eqref{norm_equivalente2} show that $\mathcal{E}$ is a topological isomorphism of $L^2(\mathbb{S}^{d-1})$ to $W(\mathbb{R}^d)$. For $\alpha,\beta>0$ we will prove that the following embedding scheme holds.
$$\begin{array}{lll}
	H^\alpha(\mathbb{S}^{d-1}) \hookrightarrow  L^2(\mathbb{S}^{d-1}) \hookrightarrow H^{-\beta}(\mathbb{S}^{d-1})  \vspace{0.3cm} \\
	\hspace{0.5cm} \downarrow \; \mathcal{E} \hspace{1.2cm} \downarrow \; \mathcal{E} \hspace{1.3cm} \downarrow  \; \mathcal{E} \hspace{1.2cm} \vspace{0.3cm} \\
	\mathcal{W}^\alpha(\mathbb{R}^{d})  \hspace{0.2cm}\hookrightarrow W(\mathbb{R}^{d}) \hspace{0.2cm} \hookrightarrow \mathcal{W}^{-\beta}(\mathbb{R}^{d}) 
	\end{array}
$$

We present first the characterization of the $\alpha$-Herglotz wave functions in terms of certain norm related to \eqref{salvador}. 
For convenience, we introduce a norm equivalent to the one given in \eqref{salvador}.
In order to do that, we notice that
using the Green-Beltrami identity (see \cite[Proposition 3.33]{A+H}), we have that for any $f\in C^2(\mathbb{S}^{d-1})$,
\begin{equation*}\label{equivalencia_tant_lapl}
\left\| \nabla_S f\right\|^2_{L^2\left( \mathbb{S}^{d-1} \right)}
=
\int_{\mathbb{S}^{d-1}} f(\theta)\left(-\Delta_S\right)f(\theta)d\sigma(\theta)
=
\left\| \left(-\Delta_S\right)^{1/2} f \right\|^2_{L^2\left( \mathbb{S}^{d-1} \right)} 
, 
\end{equation*}  
where $\Delta_S$ is the Laplace-Beltrami operator on the sphere, and 
$\left(-\Delta_S\right)^{1/2}$ is defined using the functional calculus for the positive operator $-\Delta_S$
(for a more detailed explanation see \eqref{operadores} ahead).
And therefore, the condition given in \eqref{salvador} is equivalent to
\begin{equation}
\label{salvador2}
\int_{\mathbb{R}^d}\left|u(x)+\left(-\Delta_S\right)^{1/2} u(x)\right|^2\langle x\rangle^{-3}\,dx<\infty.
\end{equation}
Here we have extended the definition of the operator $\left(-\Delta_S\right)^{1/2}$ to functions $u(x)$ defined for $x\in\mathbb{R}^d$ by using polar coordinates to write $u(x)=u(r\theta)$ with $(r,\theta)\in(0,\infty)\times\mathbb{S}^{d-1}$ and making the operator act on the spherical variable $\theta$.

This observation suggests us to introduce for any $\alpha\in\mathbb{R}$ the operators
$$\mathcal{L}^\alpha:C^{\infty}(\mathbb{S}^{d-1})^*\rightarrow C^{\infty}(\mathbb{S}^{d-1})^*,$$ 
defined, using again the the functional calculus for the positive operator $-\Delta_S$ and transposition (for a more detailed explanation see \eqref{operadores} ahead), by
\begin{equation}\label{operador_esferico}
\mathcal{L}^\alpha:=\left(I+(-\Delta_S)^{1/2}\right)^\alpha.
\end{equation} 
As before we extend the definition of $\mathcal{L}^\alpha$ to smooth functions $u(x)$ with $x\in\mathbb{R}^d$, by writing $x=r\theta$ and making $\mathcal{L}^\alpha$ acts on the spherical variable $\theta$. 
And thus, for such smooth functions $u$ and any $\alpha\in \mathbb{R}$ we can define the norm
\begin{equation}
\label{normaalpha}
 \|u\|_{\alpha}^2
	:=
	\int_{\mathbb{R}^d}  \left| \mathcal{L}^{\alpha+1}u(x)   \right|^2  \frac{dx}{<|x|>^3}.
 \end{equation} 

With this notation, from \eqref{salvador2}, we have that
the norm given in \eqref{salvador} is equivalent to $ \|\cdot\|_0$. Moreover,
\begin{equation}
\label{alpha0}
\|u\|_{\alpha}=\|\mathcal{L}^\alpha u\|_0.
\end{equation}

Our first characterization of the $\alpha$-Herglotz wave functions, given in the following theorem, is written in terms of the operators $\mathcal{L}^\alpha$ and the norms $\|\cdot\|_{\alpha}$.

\begin{theorem}\label{TEOREMA} Let $\alpha\in\mathbb{R}$. 
	An entire solution of the homogeneous Helmholtz equation $u$ is a $\alpha$-Herglotz wave function if and only if $\mathcal{L}^{\alpha}u$ is a Herglotz wave function, that is,
	\begin{equation}\label{salvador_condition}
	\|u\|_{\alpha}^2
	=
	\int_{\mathbb{R}^d}  \left| \mathcal{L}^{\alpha+1}u(x)   \right|^2  \frac{dx}{<|x|>^3}
	<\infty,
	\end{equation}
	where $\mathcal{L}^\alpha$ is defined in \eqref{operador_esferico}. Moreover, writing $u=\mathcal{E}\phi$ with $\phi\in H^\alpha(\mathbb{S}^{d-1})$, we have that
	$$\Vert u\Vert_{\alpha}\sim\Vert\phi\Vert_{ H^\alpha(\mathbb{S}^{d-1})}.$$
Hence the operator $\mathcal{E}$ defined in \eqref{definicionExtension} is a topological isomorphism of $ H^\alpha(\mathbb{S}^{d-1})$ onto $\mathcal{W}^\alpha(\mathbb{R}^d)$ provided with the norm $\|\cdot\|_{\alpha}$. 
\end{theorem}

\begin{remark}
Observe that in the case $\alpha=0$, \eqref{salvador_condition} is equivalent to \eqref{salvador} and therefore, Theorem \ref{TEOREMA} generalizes \cite[Theorem 4]{P-E+V-D}. 
\end{remark}

For the characterization of the $\alpha$-Herglotz wave functions related to the one given by Hartman-Wilcox in \eqref{HartmanWilcox} for Herglotz wave functions,  we distinguish two cases, $\alpha>0$ and $\alpha<0$. In both cases, some definitions are needed.

In case $\alpha>0$, we will use a non differential description of the Sobolev spaces in the sphere which is presented in \cite{BLP}. 

The case $0<\alpha<2$ requires the following definition.
\begin{definition}\label{def:square} Let $0<\alpha<2$. Given an integrable function $f$ on $\mathbb{S}^{d-1}$, we define the square function 
	\begin{equation*}\label{def:Sfunction}
	S_\alpha(f)^2(\theta):=\int_0^\pi\left|\frac{A_tf(\theta)-f(\theta)}{t^\alpha}\right|^2\frac{dt}{t},
	\end{equation*}
	where
	\begin{equation}\label{def:At}
	A_tf(\theta):=\frac{1}{|C(\theta,t)|}\int_{C(\theta,t)}f(\tau)\,d\sigma(\tau)
	\end{equation}
	denotes the mean of $f$ on the spherical cap centred at $\theta\in\mathbb{S}^{d-1}$ and with angle $t\in(0,\pi]$; that is,
	\begin{equation*}\label{def:casquete}
	C(\theta,t):=\{\eta\in\mathbb{S}^{d-1}:\theta\cdot\eta\geq\cos t\}.
	\end{equation*}
\end{definition}

As we did before, we extend this definition to functions $u(x)$ defined for $x\in\mathbb{R}^d$ by using polar coordinates to write $u(x)=u(r\theta)$ and making the operators $S_\alpha$ and $A_t$ act on the spherical variable $\theta$.


The next theorem gives our second characterization of the $\alpha$-Herglotz wave functions in case $0<\alpha<2$, and it is written in terms of the previous squared functions and the norm $\|\cdot\|_A$ given in \eqref{HartmanWilcox}.

\begin{theorem}\label{Alpha+1} Let $0<\alpha<2$. An entire solution of the homogeneous Helmholtz equation $u$ is a $\alpha$-Herglotz wave function if and only if $S_{\alpha}(u)$ is a Herglotz wave function, that is,
\begin{equation*}
\|S_\alpha (u)\|_{A}^{2}
=\sup_{R>0}\,\frac{1}{R}\int_{|x|<R}S_\alpha (u)^2(x)\,dx
<\infty.
\end{equation*}
Moreover, writing $u=\mathcal{E}\phi$ with $\phi\in H^\alpha(\mathbb{S}^{d-1})$, we have that
$$\|S_\alpha (u)\|_{A}\sim\Vert\phi\Vert_{H^\alpha(\mathbb{S}^{d-1})}.$$
\end{theorem}

The case $\alpha\geq2$ requires the following extension of the square function given in Definition \ref{def:square}.
\begin{definition}\label{def:square_extension} Let $\alpha\geq 2$ and $n\in\mathbb{N}$ such that  $2n\leq\alpha<2(n+1)$. Given $f,g_1,\ldots,g_n$ integrable functions on $\mathbb{S}^{d-1}$, for any $\theta\in\mathbb{S}^{d-1}$, in case $2n<\alpha<2(n+1)$, we define the square function
	\begin{equation}
	\label{def:Sgeneral}
	S_\alpha(f,g_1,\ldots,g_n)^2(\theta):=\int_0^\pi\bigg|A_tf(\theta)-f(\theta)-\sum_{k=1}^ng_k(\theta)c_k(\theta) \bigg|^2 \frac{dt }{t^{2 \alpha +1}},
	\end{equation}
	where
	$A_t$ is the mean function defined in \eqref{def:At}
	and $c_k(\theta)=A_t \big( \left|  \theta - \cdot \right|^{2k}    \big)(\theta)$, for $k=1,\ldots,n$. 
	
	And in case $\alpha=2n$, we define the square function
	\begin{align}
	S_{2n}(f,g_1,\ldots,g_n)^2(\theta):=\int_0^\pi\bigg|&A_tf(\theta)-f(\theta)-\sum_{k=1}^{n-1}\,g_k(\theta)c_k(\theta)\nonumber\\
	&-A_t g_n(\theta)c_n(\theta)   \bigg|^2 \frac{dt }{t^{4n +1}}.\label{def:Spar}
	\end{align}
\end{definition}

The next theorem gives the analogous characterization of the $\alpha$-Herglotz wave functions given in Theorem \ref{Alpha+1} for $0<\alpha<2$, but for the case $\alpha\ge 2$.
\begin{theorem}\label{Alpha+2} Let $\alpha\geq 2$ and $n\in\mathbb{N}$ such that $2n\leq \alpha<2(n+1)$. An entire solution of the homogeneous Helmholtz equation $u$ is a $\alpha$-Herglotz wave function if and only if there exist $n$ Herglotz wave functions $v_1$, $v_2, \ldots, v_n$ such that $S_{\alpha}(u,v_1,\ldots v_n)$ is a Herglotz wave function, that is,
\begin{equation}
\label{Aalpha2}
\|S_\alpha(u,v_1,\ldots v_n)\|_{A}^{2}=\sup_{R>0}\,\frac{1}{R}
\int_{|x|<R}S_\alpha(u,v_1,\ldots v_n)^2(x)|^2\,dx
<\infty.
\end{equation}
Moreover, writing $u=\mathcal{E}\phi$ with $\phi\in H^\alpha(\mathbb{S}^{d-1})$, we have that
 \begin{equation}\label{equivaNormas}
\|S_\alpha (u,v_1,\ldots v_n)\|_{A}\sim_n\Vert\phi\Vert_{H^\alpha(\mathbb{S}^{d-1})}.
	\end{equation}

\end{theorem}

In case $\alpha<0$, our characterization of the $\alpha$-Herglotz wave functions is not  written in terms of the square functions introduced above but in terms of the 
following operators.   

For convenience, we will write our characterization for $(-\alpha)$-Herglotz wave functions with $\alpha>0$ instead of $\alpha$-Herglotz wave functions with $\alpha<0$.

For $\alpha>0$, let $\mathcal{K}_\alpha $ be the integral operator defined by
\begin{equation}
\label{kalpha}
\mathcal{K}_\alpha f(\xi) := \int_{\mathbb{S}^{d-1}}K_\alpha (\xi , \eta ) f(\eta) d \sigma (\eta), \hspace{0.4cm} f \in C^\infty (\mathbb{S}^{d-1}),\ \xi\in\mathbb{S}^{d-1},
\end{equation}
where
\begin{equation*}
K_\alpha(\xi,\eta):=\frac{1}{\vert \mathbb{S}^{d-1}\vert}\int_0^1 \frac{(1-r)^{\alpha}(1+r)  }{\vert r\xi-\eta\vert^d}dr,\quad \text{for } \xi\ne\eta.
 \end{equation*}
Observe that
$$\calK_\alpha: C^{\infty}(\mathbb{S}^{d-1})\rightarrow  C^{\infty}(\mathbb{S}^{d-1}).$$

As we did before, the operators $\mathcal{K}_\alpha$ can be extended to functions 
$u(x)$ defined for $x\in\mathbb{R}^d$ by using polar coordinates to write $u(x)=u(r\theta)$ and making the operators act on the spherical variable $\theta$.

\begin{theorem} \label{alphaNegativo}
Let $\alpha>0$ and $u$ be an entire solution of the Helmholtz equation. Then $u$ is a $(-\alpha)$-Herglotz wave function if and only if $\calK_\alpha u$ is a Herglotz wave function, that is,
\begin{equation} 
\label{condicionKA}
\Vert \calK_\alpha u\Vert_A^2=\sup_{R>0}\frac{1}{R}\int_{\vert x\vert<R}\left|\calK_\alpha u(x)\right|^2 dx<\infty.
\end{equation}
Moreover, writing $u=\mathcal{E}\phi$ with 
$\phi\in H^{-\alpha}(\mathbb{S}^{d-1})$, we have that
	\begin{equation*}\label{ventana_1}
	\|\calK_\alpha u\|_A\sim  \|\phi\|_{H^{-\alpha}(\mathbb{S}^{d-1})}.
	\end{equation*}
\end{theorem}

The following proposition gives us a necessary condition to be a $(-\alpha)$-Herglotz wave function which is closer in form to \eqref{HartmanWilcox} than \eqref{condicionKA}.
\begin{proposition} \label{donsalvador_1}
	Let $\alpha\geq 0$. If $u$ is a $(-\alpha)$-Herglotz wave function, then it satisfies the condition
\begin{equation}\label{NuevaHW}
\sup_{R>0}\frac{1}{R}\int_{|x|<R}\frac{|u(x)|^2 }{(1+|x|^2)^\alpha} dx< \infty .
\end{equation} 
\end{proposition} 

Although the characterization in terms of the condition \eqref{NuevaHW} is not attained, the following partial result is obtained. 
\begin{proposition}\label{prop9}Let $\alpha>0$. If $u$ is an entire solution of the homogeneous Helmholtz equation satisfying condition \eqref{NuevaHW}, then there exists $\phi\in H^{-\beta}(\mathbb{S}^{d-1})$, with $\beta>\alpha+1/2$, such that $u=\mathcal{E}\phi.$
\end{proposition}

A result similar to Proposition \ref{donsalvador_1} but for $\alpha$-Herglotz wave functions with $\alpha>0$ is not possible. See Remark \ref{observa} below for more details. However as a consequence of Proposition \ref{donsalvador_1} and Proposition \ref{prop9}, we characterize all the entire solutions of the Helmholtz equation that are the extension operator of a distribution in the sphere.
\begin{corolario}\label{ultimo}
An entire solution $u$ of the homogeneous Helmholtz equation  is the extension operator of a distribution of the sphere if and only if 
there exists $\alpha\geq 0$ such that
\begin{equation*}
\sup_{R>0}\frac{1}{R}\int_{|x|<R}\frac{|u(x)|^2 }{(1+|x|^2)^\alpha} dx< \infty .
\end{equation*} 
\end{corolario}

We finish this section with some preliminary notions and notation that will be used through out the paper. The second section of this paper is devoted to the proofs of the results stated in this introduction, and the last section is an appendix containing several estimates of the Bessel functions. 
\subsection{Preliminary notions and notation} 
\label{preliminares}
Let $L^2(\mathbb{S}^{d-1})$ denote the space of square integrable  functions in the sphere  provided with the Lebesgue surface measure  $d\sigma$. 
It is well known (see \cite{A+H}, Chapter 2) that an orthonormal basis of this space is the set of all real valued spherical harmonics given by
$$\{Y_\ell^j /\,\ell=0,1,\ldots\ \text{and }1\leq j\leq \nu(\ell)\},$$
where 
$$
\nu (\ell)=
\frac{(2\ell+d-2)(\ell+d-3)!}{\ell!(d-2)!}.
$$

A linear operator $\mathcal{M}$  acting on formal series of spherical harmonics is call a zonal multiplier if and only if there exists a sequence of complex numbers $\{\lambda_\ell\}_{\ell=0}^\infty$ such that 
\begin{equation*}\label{multiplicadorzonal}
\mathcal{M}\left(\sum_{\ell=0}^\infty Y_\ell\right)=\sum_{\ell=0}^\infty\lambda_\ell Y_\ell,
\end{equation*}
where $Y_\ell$ is any spherical harmonics of degree $\ell$.
 We will apply this definition  to the space of distributions in the sphere, that is the dual space of $C^{\infty}(\mathbb{S}^{d-1})$, which we denote by $C^{\infty}(\mathbb{S}^{d-1})^*$. In order to do that, notice that (see \cite{Lemoine} ) every $f\in C^{\infty}(\mathbb{S}^{d-1})^*$ has a unique representation
\begin{equation}\label{representationFunct}
f=\sum_{\ell=0}^\infty\sum_{j=1}^{\nu(\ell)} \widehat{f}_{\ell j}\,Y_\ell^j, \quad \widehat{f}_{\ell j}:=\langle f,Y_\ell^j \rangle,
\end{equation}
where $\langle\cdot,\cdot\rangle$ denotes the duality in $C^\infty (\mathbb{S}^{d-1})$,
with convergence in the weak* topology and where 
\begin{equation}\label{growth}
\vert\widehat{f}_{\ell j}\vert \le C\ell^N,
\end{equation}

for some $C,N$ depending on $f$;
and conversely, every series \eqref{representationFunct} satisfying \eqref{growth} defines a distribution in the sphere. 

Moreover, when   $f\in C^{\infty}(\mathbb{S}^{d-1})$, $\widehat{f}_{\ell j}$ has the fast decay
\begin{equation*}\label{fastgrowth}
\ell^M\vert\widehat{f}_{\ell j}\vert \le C_M,
\end{equation*}
for every $M>0$, and the convergence of the series \eqref{representationFunct} is uniform in $\mathbb{S}^{d-1}$.

Hence every sequence $\lbrace \lambda_\ell\rbrace_{\ell=0}^\infty$ with polynomial growth, that is,
\begin{equation*}
\vert\lambda_{\ell}\vert \le C\ell^N,
\end{equation*}
for some $N$, defines a zonal multiplier from  $C^\infty (\mathbb{S}^{d-1})$ into  $C^\infty (\mathbb{S}^{d-1})$
and by transxposition, also from $C^\infty (\mathbb{S}^{d-1})^*$ into $C^\infty (\mathbb{S}^{d-1})^*$.

Notice that the composition of two zonal multiplier defined on $C^{\infty}(\mathbb{S}^{d-1})^*$ is another zonal multiplier, and the composition commutes.

The space of all spherical harmonics of degree $\ell$ spanned by  $\{Y_\ell^j :1\leq j\leq \nu(\ell)\}$ is the eigenspace of $-\Delta_S$ corresponding to the eigenvalue $\ell(\ell+d-2)$. 
Thus, $-\Delta_S$ is a zonal multiplier with associated sequence $\{\ell(\ell+d-2)\}_{\ell=0}^\infty.$

From here, we can see the operators $(-\Delta_S)^{1/2}$ and $\mathcal{L}^\alpha$ introduced in \eqref{salvador2} and \eqref{operador_esferico} as zonal multipliers defined on $C^{\infty}(\mathbb{S}^{d-1})^*$ with 
associated sequences 
\begin{equation}
\label{operadores}
\left\{\ell^{1/2}(\ell+d-2)^{1/2}\right\}_{\ell=0}^\infty
\qquad
\text{and}
\qquad
\left\{
\left(1+\ell^{1/2}(\ell+d-2)^{1/2}\right)^\alpha
\right\}_{\ell=0}^\infty
\end{equation}
respectively.

Moreover, for $\alpha\in\R,$ we can define the Sobolev space  $H^\alpha (\mathbb{S}^{d-1})$
as  the space of all distributions $f\in C^{\infty}(\mathbb{S}^{d-1})^*$ such that $\mathcal{L}^\alpha f\in  L^2(\mathbb{S}^{d-1})$, and equipped with the norm 
\begin{equation}\label{def:H}
\|f\|_{H^{\alpha}(\mathbb{S}^{d-1})}^2
=
\sum_{\ell=0}^\infty\sum_{j=1}^{\nu(\ell)} \left(1+\ell^{1/2}(\ell+d-2)^{1/2}\right)^{2\alpha}
|\widehat f_{\ell j}|^2
<\infty.
\end{equation}

On the other hand, since any $f\in C^{\infty}(\mathbb{S}^{d-1})^*$ admits representation \eqref{representationFunct}, using the continuity of the Fourier transform in temperate distributions and the Funk-Hecke's formula (see \cite[pp.37]{Helgason} and also \cite[Lemma 4]{P-E+V-D}), and writing $x=r\theta$,  we have that 
\begin{align}
\mathcal{E}f(x)=&\sum_{\ell=0}^\infty\sum_{j=1}^{\nu(\ell)} \widehat{f}_{\ell j}\,\mathcal{E}Y_\ell^j(r\theta)\nonumber\\
=&\sum_{\ell=0}^\infty\sum_{j=1}^{\nu(\ell)} (2\pi)^{\frac{1}{2}} i^{\ell} 
\frac{J_{\mu(\ell)}(r)}{r^{\frac{d-2}{2}}}
\widehat{f}_{\ell j}Y_\ell^j(\theta)\label{repW},
\end{align}
where $\mu (\ell)= \ell + \frac{d-2}{2}.$ And therefore, $\mathcal{E}$ acts as a zonal multiplier in the spherical variable $\theta$.

Notice that the fast uniform convergence to zero of the Bessel functions on compact sets, implies that the series \eqref{repW} converges uniformly on compact subsets of $\R^d$.

Throughout this paper, for $X,Y\geq 0$, we will write $X\sim Y$ if there exists a constant $c>0$, depending on at most the dimension $d$, such that $c^{-1}Y \leq X \leq cY$, and $X \lesssim Y$ if there exists a similar uniform constant $c$ such that $X \leq c Y$. We will write $X \sim_b Y$ if the constants above depend on a specific parameter $b>0$.



\section{Proofs}
We start this section proving Theorem \ref{TEOREMA}.

{\emph{Proof of Theorem \ref{TEOREMA}.}} 
Let $u$ be a $\alpha-$Herglotz wave function. Then there exists $\phi \in H^\alpha (\mathbb{S}^{d-1})$ such that  $u= \mathcal{E}\phi$. 
And therefore, since from \eqref{operadores} and \eqref{repW}, we have that
$ \mathcal{L}^\alpha$ and $\mathcal{E}$ are zonal multipliers, we can write
$$\mathcal{L}^\alpha u= \mathcal{L}^\alpha \mathcal{E} \phi= \mathcal{E} \mathcal{L}^\alpha  \phi.$$ 
And thus, $\mathcal{L}^\alpha u$ is a  Herglotz wave function since $\mathcal{L}^\alpha  \phi \in L^2(\mathbb{S}^{d-1})$.

Conversely, let $u$ be an entire solution of the homogeneous Helmholtz equation such that $\mathcal{L}^\alpha u$  is a Herglotz wave function. Then, there exists $\phi \in L^2(\mathbb{S}^{d-1})$ such that $ \mathcal{L}^\alpha u= \mathcal{E}\phi$, and therefore,
\begin{equation*}
u=\mathcal{L}^{-\alpha } \mathcal{E} \phi =\mathcal{E} \mathcal{L}^{-\alpha}\phi.
\end{equation*}
From here, we have that $u$ is a $\alpha-$Herglotz wave function  since $\mathcal{L}^{-\alpha}f \in H^\alpha (\mathbb{S}^{d-1})$.

On the other hand, writing $u=\mathcal{E}\phi$ with $\phi\in H^\alpha(\mathbb{S}^{d-1})$ and using \eqref{alpha0}, \eqref{norm_equivalente2}, \eqref{operadores} and \eqref{def:H}, we get
\begin{equation}
\label{equinorms}
\|u\|_\alpha=\|\mathcal{E}\phi\|_\alpha
=\|\mathcal{L}^\alpha\mathcal{E}\phi\|_0
=\|\mathcal{E}\mathcal{L}^\alpha\phi\|_0
\sim\|\mathcal{L}^\alpha\phi\|_{L^2(\mathbb{S}^{d-1})}
=\|\phi\|_{H^\alpha(\mathbb{S}^{d-1})}.
\end{equation}
Finally, since by definition $\mathcal{E}:L^2(\mathbb{S}^{d-1})\rightarrow\mathcal{W}$ and  $\mathcal{L}^{\alpha}: H^{\alpha}(\mathbb{S}^{d-1})\rightarrow L^2(\mathbb{S}^{d-1}) $ are bijections, we have that
 $$\mathcal{E}=\mathcal{L}^{-\alpha} \mathcal{E}\mathcal{L}^{\alpha}: H^\alpha(\mathbb{S}^{d-1})\rightarrow\mathcal{W}^\alpha(\mathbb{R}^d)$$ 
is also a bijection, and thus, from \eqref{equinorms}, $\mathcal{E}$ is a topological isomorphism of $H^\alpha(\mathbb{S}^{d-1})$ to $\mathcal{W}^\alpha(\mathbb{R}^d)$.
\hfill$\Box$

The proof of theorems \ref{Alpha+1} and \ref{Alpha+2}
requires the following auxiliary results that can be found in \cite{BLP}.
\begin{theorem}\label{teorema1}{\rm(\cite[Theorem 1.1 and Theorem 1.2]{BLP}).}
	Let $\alpha>0$, $n$ the non-negative integer number such that  $2n\leq\alpha<2(n+1)$, and $\phi\in L^2(\mathbb{S}^{d-1})$.
	\begin{enumerate}[label=\arabic*)]
		\item 
		If $n=0$, then 
		$\phi \in H^\alpha (\mathbb{S}^{d-1})$ 
		if an only if 
		$S_\alpha (\phi) \in L^2(\mathbb{S}^{d-1})$. Moreover,
		\[\|\phi\|_{H^{\alpha}(\mathbb{S}^{d-1})}\sim \|S_\alpha (\phi)\|_{L^2(\mathbb{S}^{d-1})}.\]
		\item 
		If $n\geq 1$, then 
		$\phi \in H^\alpha \left(  \mathbb{S}^{d-1} \right)$ 
		if and only if there exist  
		$\phi_1,\phi_2, \cdot \cdot \cdot, \phi_n \in L^2 \left( \mathbb{S}^{d-1}  \right) $ 
		such that  
		$S_\alpha (\phi,\phi_1, \phi_2, \cdot \cdot \cdot, \phi_n )  \in  L^2(\mathbb{S}^{d-1})$. Moreover,
		\begin{equation}\label{bombilla}	
		\|\phi\|_{H^{\alpha}(\mathbb{S}^{d-1})}\sim_{n}\|S_\alpha (\phi,\phi_1, \phi_2, \cdot \cdot \cdot, \phi_n )\|_{L^2(\mathbb{S}^{d-1})}.
		\end{equation}
	\end{enumerate}
\end{theorem}

\begin{lemma}\label{lemma:At}{\rm (\cite[Lemma 2.1]{BLP}).} 
	For each $t \in (0,\pi]$, the operator $A_t$ defined in \eqref{def:At} is a zonal Fourier multiplier with associated sequence $\{m_{\ell,t}\}_{\ell=0}^\infty$ given by 
\begin{equation}\label{mutliplierAt}
	m_{\ell,t}:=C_{t,d}\int_{\cos t}^1P_{\ell,d}(s)(1-s^2)^{\frac{d-3}{2}}\,ds,\;\;\; \ell= 0,1, \cdots
	\end{equation}
where $C_{t,d}=\frac{|\mathbb{S}^{d-2}|}{|C(\xi,t)|}$ and $P_{\ell,d}$ denotes the Legendre polynomial of degree $\ell$ in $d$ dimensions.
\end{lemma}
The proofs of Theorems \ref{Alpha+1} and \ref{Alpha+2} are similar. We omit the proof of the Theorem \ref{Alpha+1} because it is simpler than the proof of the Theorem \ref{Alpha+2}, since in the first case, no summation term is involved in the definition of the square function (see Definitions \ref{def:square} and \ref{def:square_extension}). 

\emph{Proof of Theorem \ref{Alpha+2}.} 
We consider first the case $2n<\alpha<2(n+1)$. 

We start with the necessary condition. Let $u$ be a $\alpha-$Herglotz wave function, then $u= \mathcal{E}\phi$ with $\phi \in H^\alpha (\mathbb{S}^{d-1})$, and thus, from Theorem \ref{teorema1}, there exist
$\phi_1,\phi_2, \ldots, \phi_n \in L^2 \left( \mathbb{S}^{d-1}  \right) $ such that
\begin{equation}\label{inter1}
\|\phi\|_{H^{\alpha}(\mathbb{S}^{d-1})}\sim_{n}\|S_\alpha (\phi,\phi_1, \phi_2, \ldots, \phi_n )\|_{L^2(\mathbb{S}^{d-1})}.
\end{equation}
We introduce the Herglotz wave functions 
\[v_k:=\mathcal{E}\phi_k,\quad k=1,\ldots,n.\] 
From \eqref{def:Sgeneral}, we have that
\begin{eqnarray}
S_\alpha(u,v_1,\ldots v_n)^2(x)
&=&\int_0^\pi\bigg|A_tu(x)-u(x)-\sum_{k=1}^nv_k(x)c_k(t)  \bigg|^2
\frac{dt}{t^{2\alpha+1}}
\nonumber\\
&=&
\int_0^\pi
\bigg|\mathcal{E}\phi_t(x)\bigg|^2
\frac{dt}{t^{2\alpha+1}},
\label{paso1}
\end{eqnarray}
where
\begin{equation}\label{density}
\phi_t:=A_t\phi-\phi-\sum_{k=1}^n\phi_kc_k(t),\quad t\in(0,\pi].
\end{equation}
Observe that in \eqref{paso1} we have used the fact that $\mathcal{E}$ and $A_t$ are zonal multipliers (see \eqref{repW} and Lemma \ref{lemma:At} respectively) and therefore, they commute.

From \eqref{paso1}, for $R>0$ fixed, we have that
\begin{equation}
\label{paso2}
\frac{1}{R}
\int_{|x|<R} S_\alpha(u,v_1,\ldots v_n)^2(x)\,dx=
\int_0^\pi\frac{1}{R}\int_{|x|<R}\bigg|\mathcal{E}\phi_t(x)  \bigg|^2dx\,\frac{dt}{t^{2\alpha+1}}.
\end{equation}
Since $\phi\in L^2(\mathbb{S}^{d-1})$, for each $t\in(0,\pi)$, $\phi_t\in L^2(\mathbb{S}^{d-1})$, and therefore, $\mathcal{E}\phi_t$ is a Herglotz wave function. Thus, using \eqref{norm_equivalente} in \eqref{paso2} we get
$$
\frac{1}{R}
\int_{|x|<R}S_\alpha(u,v_1,\ldots v_n)^2(x)\,dx\lesssim \int_0^\pi\|\phi_t\|_{L^2(\mathbb{S}^{d-1})}^2\frac{dt}{t^{2\alpha+1}}.
$$
From here, using \eqref{density}, \eqref{def:Sgeneral} and \eqref{inter1} we obtain
\begin{align*}
\frac{1}{R}
\int_{|x|<R} S_\alpha(u,v_1,\ldots v_n)^2(x)\,dx&\lesssim \|S_\alpha(\phi,\phi_1,\ldots \phi_n)\|_{L^2(\mathbb{S}^{d-1})}^2\nonumber\\
&\sim_{n}\|\phi\|_{H^\alpha(\mathbb{S}^{d-1})}^2.
\end{align*}
Then, taking the supremum in $R$, since $\phi \in H^\alpha (\mathbb{S}^{d-1})$, we deduce that 
\begin{align}\label{dir1}
\|S_\alpha(u,v_1,\ldots v_n)\|_A\lesssim_{n}\|\phi\|_{H^\alpha(\mathbb{S}^{d-1})}<\infty.
\end{align}

Now we will prove the sufficient condition. First we will assume that $u$ is a Herglotz wave function, and then we will extend the result to entire solutions of the homogeneous Helmholtz equation.

For a Hergltoz wave function $u$, using the sufficient condition, we have that there exist $f,g_1,\ldots,g_n\in L^2(\mathbb{S}^{d-1})$ such that
$$
u=\mathcal{E}f,v_1=\mathcal{E}g_1,\ldots,v_n=\mathcal{E}g_n.
$$

From \eqref{def:Sgeneral}, using Lemma \ref{lemma:At}, we can write
\begin{align*}
S_\alpha(f,g_1,\ldots,g_n)^2(\theta)
&=
\int_0^\pi \left|A_t f(\theta)-f(\theta)-\sum_{k=1}^n g_k(\theta) c_k (t)\right|^2 \frac{dt}{t^{2 \alpha +1}}
\\
&=
\int_0^\pi \left| \sum_{\ell=0}^\infty\sum_{j=1}^{\nu(\ell)}\widehat{h}_{\ell j}(t) Y_\ell^j ( \theta)\right|^2\frac{dt}{t^{2 \alpha +1}},
\end{align*}
with
\begin{equation*}\label{def:h}
\widehat{h}_{\ell j}(t):= m_{\ell,t}\widehat{f}_{\ell j}-\widehat{f}_{\ell j}-\sum_{k=1}^n\widehat{g_k}_{\ell j}c_k(t),
\end{equation*}
where $m_{\ell,t}$ is defined in \eqref{mutliplierAt},
and therefore, we have that
\begin{equation}\label{formulaS}
\|S_\alpha(f,g_1,\ldots,g_n)\|_{L^2(\mathbb{S}^{d-1})}^2=\int_0^\pi \sum_{\ell=0}^\infty\sum_{j=1}^{\nu(\ell)}|\widehat{h}_{\ell j}(t)|^2\frac{dt}{t^{2\alpha+1}}.
\end{equation}
Arguing in a similar way, but using \eqref{repW}, we can write
\begin{align*}
S_\alpha(u,v_1,\ldots,v_n)^2(x)
&=
\int_0^\pi \left|A_t u(x)-u(x)-\sum_{k=1}^n v_k(x) c_k (t)\right|^2 \frac{dt}{t^{2 \alpha +1}}
\\
&=
\int_0^\pi 2\pi\left| \sum_{\ell=0}^\infty\sum_{j=1}^{\nu(\ell)}i^\ell\,\widehat{h}_{\ell j}(t)\frac{J_{\mu(\ell)}(r)}{r^{\frac{d-2}{2}}} Y_\ell^j ( \theta)\right|^2\frac{dt}{t^{2 \alpha +1}},
\end{align*}
where $x=r\theta$, and therefore, for $R>0$ fixed, we have that
\begin{align}
\frac{1}{R} \int_{|x|<R} S_\alpha (u,v_1,\ldots, v_n)^2(x)dx&\nonumber\\
&\hspace*{-2.5cm}=\sum_{\ell=0}^\infty\sum_{j=1}^{\nu(\ell)}\frac{2\pi}{R}\int_0^R|J_{\mu(\ell)}(r)|^2r\,dr\int_0^{\pi}|\widehat{h}_{\ell j}(t)|^2\frac{dt}{t^{2\alpha+1}}.\label{segundo}
\end{align}

From the following asymptotic formula (see \cite[pp. 134]{LE})
$$J_\mu(r)=\sqrt{\frac{2}{\pi r}} \cos \left( r-\frac{\mu \pi}{2}-\frac{\pi}{4}    \right)+ O(r^{-3/2}), \hspace{0.3cm} r \rightarrow \infty,$$ 
we get (see \cite{CK} for the case $\mu=\ell+\frac{1}{2}, \; \ell \in \mathbb{N}$)
\begin{equation}
\label{limite}
\lim_{R \rightarrow \infty}\frac{1}{R} \int_0^R |J_\mu(r)|^2rdr \sim 1, \hspace{0.3cm} \mu \geq 0.
\end{equation}
Using \eqref{limite} in \eqref{segundo}, the Fatou's lemma, \eqref{formulaS} and \eqref{bombilla}, we obtain
\begin{align}
\|S_\alpha(u,v_1,\ldots v_n)\|_A^2
&\sim
\sum_{\ell=0}^\infty\sum_{j=1}^{\nu(\ell)}\int_0^{\pi}|\widehat{h}_{\ell j}(t)|^2\frac{dt}{t^{2\alpha+1}}
\nonumber
\\
&\ge
\int_0^{\pi}\sum_{\ell=0}^\infty\sum_{j=1}^{\nu(\ell)}|\widehat{h}_{\ell j}(t)|^2\frac{dt}{t^{2\alpha+1}}
\nonumber
\\
&=
\|S_\alpha(f,g_1,\ldots,g_n)\|_{L^2(\mathbb{S}^{d-1})}^2
\nonumber
\\
&\sim_{n}
\|f\|_{H^\alpha(\mathbb{S}^{d-1})}^2.
\label{tercero}
\end{align}
From \eqref{tercero} and the hypothesis \eqref{Aalpha2}, we conclude that $u=\mathcal{E}f$ with $f\in H^\alpha(\mathbb{S}^{d-1})$, that is, $u$ is a $\alpha-$Herglotz wave function.

Moreover, \eqref{tercero} together with \eqref{dir1} gives \eqref{equivaNormas}, which completes the proof whenever $u$ is a Herglotz wave function.

Now we consider $u$ to be an entire solution of the homogeneous Helmholtz equation satisfying \eqref{Aalpha2}. 

As it is well known (see \cite[Chapter 3]{CK} for $d=3$ and \cite{P-E+V-D} for greater dimensions), writing $x=r\theta$, any entire solution of the homogeneous Helmholtz equation $u$ can be expanded  as
\begin{equation}
\label{expansionu}
u(x)=(2 \pi)^{1/2}
\sum_{\ell=0}^\infty\sum_{j=1}^{\nu(\ell)}i^\ell a_{\ell j}\, \frac{J_{\mu(\ell)}(r)}{r^{\frac{d-2}{2}}}Y_\ell^j(\theta),
\qquad x\in\mathbb{R}^d,
\end{equation}
for certain coefficients $a_{\ell j}$, with uniform convergence in compact sets.

For each $N\in \N$, let $\pi_N$ be the orthogonal projection of $L^2(\mathbb{S}^{d-1})$ onto the space of spherical harmonics of degree less than or equal to $N$ given by
$$\pi_N \left( \sum_{\ell=0}^\infty\sum_{j=1}^{\nu(\ell)}\widehat{g}_{\ell j} Y_\ell^j\right)=\sum_{\ell=1}^N\sum_{j=1}^{\nu(\ell)}\widehat{g}_{\ell j} Y_\ell^j.$$
By means of \eqref{expansionu}, we can extend this definition to $u(x)$ with $x\in\mathbb{R}^d$, considering $\pi_N$ acting on the spherical variable.

Since $\pi_N u$ and $\pi_N v_1,\ldots,\pi_N v_n$  are Herglotz  wave functions, we can write
$$\pi_N u =\mathcal{E}f_N \text{ and } \pi_N v_1=\mathcal{E}\pi_N g_1,\ldots,\pi_Nv_n=\mathcal{E}\pi_N g_n,$$
where 
\begin{equation}
\label{fN}
f_N(\theta)=\sum_{\ell=0}^N\sum_{j=1}^{\nu(\ell)}a_{\ell j} Y_\ell^j(\theta)
\text{ and }
v_1=\mathcal{E}g_1,\ldots, v_n=\mathcal{E}g_n.
\end{equation}
Since we have proven the result for Herglotz wave functions, we have that
\begin{equation}
\label{paraHerglotz}
\|f_N\|_{H^{\alpha}(\mathbb{S}^{d-1})}
\lesssim_{n}
\|S_\alpha(\pi_N u,\pi_N v_1,\ldots, \pi_N v_n)\|_A.
\end{equation}

On the other hand, it is easy to check that
\begin{equation*}
\label{easy}
\|S_\alpha(\pi_N u,\pi_N v_1,\ldots \pi_N v_n)\|_A\leq \|S_\alpha(u,v_1,\ldots ,v_n)\|_A,
\end{equation*}
and using this inequality in \eqref{paraHerglotz}, from the hypothesis \eqref{Aalpha2},  we get
\begin{equation*}
\label{parau}
\|f_N\|_{H^{\alpha}(\mathbb{S}^{d-1})}
\lesssim_{n}
\|S_\alpha(u,v_1,\ldots ,v_n)\|_A<\infty.
\end{equation*}
From here, taking into account \eqref{fN} and \eqref{expansionu}, we conclude that
$$
f:=\lim_{N\rightarrow\infty} f_N=\sum_{\ell=1}^\infty\sum_{j=1}^{\nu(\ell)}a_{\ell j} Y_\ell^j\in H^{\alpha}(\mathbb{S}^{d-1}),
$$
$u=\mathcal{E}f$ 
and $\Vert f\Vert_{H^{\alpha}(\mathbb{S}^{d-1})}\lesssim_{n}  \|S_\alpha(u,v_1,\ldots, v_n)\|_A$.
Hence the proof is also completed for $u$ being an entire solution of the Helmholtz equation.

The proof in the case $\alpha=2n$ is similar to the previous one, but replacing $S_\alpha(u,v_1,\ldots, v_n)$ by $S_{2n}(u,v_1,\ldots, v_n)$, the square function defined in \eqref{def:Spar}.
\begin{flushright}
$\Box$
\end{flushright}

In order to prove Theorem \ref{alphaNegativo} we consider $\alpha>0$ and introduce the zonal multipliers 
$\mathcal{M}_\alpha$ defined on $C^\infty (\mathbb{S}^{d-1})^*$ by the sequence $\{ \beta(\alpha, \ell +1)  \}_{\ell =0}^\infty$, where
$$ \beta(x,y)=\int_0^1 t^{x-1}(1-t)^{y-1 }dt, \; \; \; \; x>0,\;  y>0,$$ is the beta function.
$\mathcal{M}_\alpha$ can be extended to functions $u$ defined on $\mathbb{R}^d$. 
Writing the beta function in terms of the gamma function and using the Striling's formula, one can see that
\begin{equation}\label{asintoticabeta}
\beta (\alpha , \ell +1 ) \sim_\alpha \ell^{-\alpha}, \qquad \ell \rightarrow \infty,
\end{equation} 
and therefore $\mathcal{M}_\alpha^{-1}$ is also a zonal Fourier multiplier on $C^\infty (\mathbb{S}^{d-1})^*$.

On the other hand, from \eqref{def:H} we see that a distribution $f\in C^\infty(\mathbb{S}^{d-1})^*$, which admits the representation \eqref{representationFunct}, belongs to  $H^{-\alpha }(\mathbb{S}^{d-1})$ if and only if 
$$\sum_{\ell=0}^\infty\sum_{j=1}^{\nu(\ell)} \widehat{f}_{\ell j}(1+\ell)^{-\alpha} Y_\ell^j\in L^2(\mathbb{S}^{d-1}),$$
and taking into account \eqref{asintoticabeta}, this holds if and only if $\mathcal{M}_\alpha f \in L^2(\mathbb{S}^{d-1})$. Moreover,
\begin{equation}\label{ventana}
\|f\|_{H^{-\alpha }(\mathbb{S}^{d-1})} \sim \| \mathcal{M}_\alpha f\|_{L^{2 }(\mathbb{S}^{d-1})}.
\end{equation}

The following lemma relates the zonal multiplier $\mathcal{M}_\alpha$ and the integral operator $\mathcal{K}_\alpha$ introduced in \eqref{kalpha}, and we will use it in the proof of Theorem \ref{alphaNegativo}.

\begin{lemma} \label{silla}
Let be $\alpha>0,$
$\mathcal{K}_\alpha$ the integral operator defined for any $\phi\in C^\infty (\mathbb{S}^{d-1})$ by \eqref{kalpha}, and 
$\mathcal{M}_\alpha$ the zonal multiplier defined by the sequence $\{ \beta(\alpha, \ell +1)  \}_{\ell =0}^\infty$, where
$\beta(x,y)$ denotes the beta function.
It holds that
\begin{equation}\label{silla_1}
\mathcal{K}_\alpha \phi =\mathcal{M}_\alpha \phi, \qquad \phi \in C^\infty (\mathbb{S}^{d-1}).
\end{equation}
\end{lemma}
Before proving this lemma, we present the proof of Theorem \ref{alphaNegativo}.

\emph{Proof of Theorem \ref{alphaNegativo}.}  Let $u$ be a $(-\alpha)$-Herglotz wave function then, there exists $\phi \in H^{-\alpha} (\mathbb{S}^{d-1})$ such that $u=\mathcal{E}\phi$.   

Since  $\mathcal{M}_\alpha$  and $\mathcal{E}$ are Fourier multipliers (see \eqref{repW}), they commute, and therefore, using  (\ref{silla_1}), we can write 
\begin{equation}
\label{kem}
\mathcal{K}_\alpha u= \mathcal{M}_\alpha u= \mathcal{M}_\alpha \mathcal{E} \phi=\mathcal{E} \mathcal{M}_\alpha \phi.
\end{equation}
Notice that the first identity in the previous expression holds because $u=\mathcal{E}\phi$ is a $C^\infty$ function in $\mathbb{R}^d$.

Since $\phi\in H^{-\alpha}(\mathbb{S}^{d-1})$, from \eqref{ventana}, $\mathcal{M}_\alpha \phi \in L^2(\mathbb{S}^{d-1})$, and therefore, from \eqref{kem}, 
$ \mathcal{K}_\alpha u $ is a Herglotz wave function.
 
Conversely, if $\mathcal{K}_\alpha u$ is a Herglotz wave function, there exists $\psi \in L^2 (\mathbb{S}^{d-1})$ such that $\mathcal{K}_\alpha u=\mathcal{E} \psi$. From \eqref{ventana}, $\phi= \mathcal{M}_\alpha^{-1} \psi
\in H^{-\alpha}(\mathbb{S}^{d-1})$, and we can write
\begin{equation}
\label{reciproco}
\mathcal{M}_\alpha u=\calK_\alpha u=\mathcal{E} \psi=\mathcal{E}\mathcal{M}_\alpha \phi=\mathcal{M}_\alpha \mathcal{E} \phi,
\end{equation}
and thus
$u=\mathcal{E} \phi$ is a $(-\alpha)$-Herglotz wave function.

Notice that the first identity in \eqref{reciproco} holds because by hypothesis, $u$ is an entire solution of the homogeneous Helmholtz equation, and therefore, is a $C^\infty$ function in $\mathbb{R}^d$.

Finally, using \eqref{ventana}, \eqref{norm_equivalente}, the fact that $u=\mathcal{E} \phi$ is a $C^\infty$ function in $\mathbb{R}^d$ and \eqref{silla_1}, we get
$$\|\phi \|_{H^{-\alpha}(\mathbb{S}^{d-1})} 
\sim \|\mathcal{M}_\alpha \phi \|_{L^{2} (\mathbb{S}^{d-1})}  
\sim \| \mathcal{E} \mathcal{M}_\alpha \phi \|_{A} 
= \| \mathcal{M}_\alpha u \|_A 
=\| \mathcal{K}_\alpha u \|_A.
$$
\hfill$\Box$

It remains to prove Lemma \ref{silla}.
In order to do this, we need some properties of the Poisson transform, defined for $0\le r<1$ by
\begin{equation}\label{Poisson}
P_rf(\xi)=\int_{\mathbb{S}^{d-1}}p_r(\xi , \eta) f(\eta) d \sigma (\eta), \qquad f \in L^1(\mathbb{S}^{d-1}),\xi\in \mathbb{S}^{d-1},
\end{equation}
where
$$p_r(\xi, \eta)=\frac{1}{|\mathbb{S}^{d-1}|}\ \frac{1-r^2}{|r \xi -\eta  |^d}, 
\qquad\xi, \eta \in \mathbb{S}^{d-1}.$$
It is well known (see \cite{ABR}) that 
$P_rf \in C^{\infty}(\mathbb{S}^{d-1})$, and
$$p_r(\xi, \eta)=\sum_{\ell =0}^\infty \sum_{j=1}^{\nu (\ell)}r^\ell Y_\ell^j (\xi) Y_\ell^j (\eta),$$
which implies that $P_r$ is a zonal multiplier with associated sequence $\{ r^\ell  \}_{\ell =0}^\infty$.

\emph{Proof of Lemma \ref{silla}.}  
Given $\phi\in C^\infty(\mathbb{S}^{d-1})$, we can write
$$ 
\phi=\sum_{\ell=0}^\infty\sum_{j=1}^{\nu(\ell)} \widehat{\phi}_{\ell j}\,Y_\ell^j.
$$
and since $P_r$ is a zonal multiplier with associated sequence $\{ r^\ell  \}_{\ell =0}^\infty$, we have that
$$P_r \phi=\sum_{\ell=0}^\infty\sum_{j=1}^{\nu(\ell)} r^{\ell}\widehat{\phi}_{\ell j}\,Y_\ell^j,$$
From here, using \eqref{Poisson} and the fast decay of $\widehat{\phi}_{\ell j}$, for any $\xi\in\mathbb{S}^{d-1}$, we obtain that
\begin{equation}
\label{uno}
 \int_0^1 (1-r)^{\alpha-1}P_r \phi(\xi) dr=\sum_{\ell=0}^\infty\sum_{j=1}^{\nu(\ell)} \beta(\alpha,\ell+1)\widehat{\phi}_{\ell j}\,Y_\ell^j(\xi)=\mathcal{M}_\alpha\phi(\xi).
\end{equation}  
On the other hand,
using \eqref{Poisson}, Fubini's theorem and \eqref{kalpha}, for any $\xi\in\mathbb{S}^{d-1}$, we can write
\begin{equation}
\label{dos}
\int_0^1 (1-r)^{\alpha-1}P_r \phi(\xi)dr=\int_{\mathbb{S}^{d-1}}K_\alpha(\xi,\eta)\phi(\eta)d\sigma(\eta)=\calK_\alpha \phi(\xi). 
\end{equation}
From \eqref{uno} and \eqref{dos}, identity \eqref{silla_1} holds.  \hfill$\Box$

The proof of Proposition \ref{donsalvador_1} requires the following lemma involving Bessel functions, that will be proven in the Appendix.

\begin{lemma} \label{lema_principal}
	Let $\mu \geq 1/2 $ and $s \geq 0$. It holds
	\begin{equation}\label{equivalencia_condicion_Herglotz}
	\sup_{R>0} \frac{1}{R}\int_{0}^R   \left|  J_\mu (r)  \right|^2 \frac{r\,dr}{(1+r^2)^s}  \sim
	(1+\mu)^{2-2s}  \int_0^\infty     \left|   J_\mu (r) \right|^2 \frac{r\,dr}{<r>^3},
	\end{equation}
where $\langle r\rangle=(1+r^2)^{1/2}$.	
\end{lemma}

\emph{Proof of Proposition \ref{donsalvador_1}.} 
Let $u$ be a $(-\alpha)$-Herglotz wave function. Then, there exists $\phi\in H^{-\alpha}(\mathbb{S}^{d-1})$ such that
$u=\mathcal{E}\phi$, and from (\ref{repW}) we can write
\begin{equation}
\label{u}
u(x)=(2\pi)^{\frac{1}{2}}
\sum_{\ell=0}^\infty\sum_{j=1}^{\nu(\ell)}i^\ell \widehat{\phi}_{\ell j}\,
\frac{J_{\mu(\ell)}(r)}{r^{\frac{d-2}{2}}}\, Y_\ell^j(\theta),
\qquad
x=r\theta.
\end{equation}
From \eqref{normaalpha}, using polar coordinates, \eqref{u}, and \eqref{operadores}, we get
\begin{align*}
\|u\|_{-\alpha}^2
&= \int_{\mathbb{R}^d}\left|\mathcal{L}^{-\alpha+1}u(x)   \right|^2\frac{dx}{<|x|>^3}\\
&= 
\int_0^\infty  \sum_{\ell=0}^\infty
\sum_{j=1}^{\nu(\ell)} 
\left|\widehat{\phi}_{\ell j}\right|^2 
\left( 1+\ell^{1/2}(\ell+d-2)^{1/2}\right)^{2- 2 \alpha} 
\left| J_{\mu (\ell)}(r) \right|^2 
\frac{r\,dr}{<r>^3}.
\end{align*}	
Since $u$ is a $(-\alpha)$-Herglotz wave function, from Theorem \ref{TEOREMA} and the dominated convergence theorem, we have that 
$$
\sum_{\ell=0}^\infty\sum_{j=1}^{\nu(\ell)} 
\left|\widehat{\phi}_{\ell j}\right|^2 
\left( 1+\ell^{1/2}(\ell+d-2)^{1/2}\right)^{2- 2 \alpha} 
\int_0^\infty 
\left| J_{\mu (\ell)}(r) \right|^2 
\frac{r\,dr}{<r>^3}
\le C,
$$
where $C$ is an absolute constant.
In particular, for any $N \in \mathbb{N}$, it holds
$$
\sum_{\ell=0}^N\sum_{j=1}^{\nu(\ell)} 
\left|\widehat{\phi}_{\ell j}\right|^2 
\left( 1+\ell^{1/2}(\ell+d-2)^{1/2}\right)^{2- 2 \alpha} 
\int_0^\infty 
\left| J_{\mu (\ell)}(r) \right|^2 
\frac{r\,dr}{<r>^3}
\le C,
$$
From here, using \eqref{equivalencia_condicion_Herglotz} with $s=\alpha\ge0$, since $\mu(\ell)=\ell+\frac{d-2}{2}$, for any $R>0$ we obtain
$$
\sum_{\ell=0}^N \sum_{j=1}^{\nu(\ell)}
\left|\widehat{\phi}_{\ell j}\right|^2 
\frac{1}{R}\int_0^R  \left| J_{\mu (\ell)}(r) \right|^2 \frac{r\,dr}{(1+r^2)^\alpha}\leq C.
$$
And therefore, using again the dominated convergence theorem, from \eqref{u},we have that
$$
\frac{1}{R}\int_0^R
\sum_{\ell=0}^\infty\sum_{j=1}^{\nu(\ell)}
\left|\widehat{\phi}_{\ell j}\right|^2 
\left| J_{\mu (\ell)}(r) \right|^2
\frac{r\,dr}{(1+r^2)^\alpha} =\frac{1}{R}\int_{|x|<R}\frac{|u(x)|^2}{(1+|x|^2)^\alpha}\,dx
<\infty.
 $$ 
 \begin{flushright}
$\Box$
 \end{flushright}			

The following remark shows that a result similar to Proposition \ref{donsalvador_1} is not true for $\alpha$-Herglotz wave functions with $\alpha>0.$
\begin{remark}\label{observa} 
Let $\alpha >0$ and $u$ be a $\alpha$-Herglotz wave function. Then $u$ can not satisfy the growth condition
\begin{equation}\label{pincel}
\sup_{R>0} \frac{1}{R} \int_{|x|<R}|u(x)|^2(1+|x|^2)^\alpha dx < \infty.
\end{equation}
This can be proved by a contradiction argument. Assume that $u$ satisfies \eqref{pincel} then, for any $R>0$ we have that
$$
\frac{1}{R}\int_{R/2 < |x| <R} |u(x)|^2\left(  1+|x|^2  \right)^\alpha dx \leq C,
$$ 
where $C$ is an absolute constant independent of $R,$ and therefore, since $\alpha>0,$ we get
\begin{equation}
\label{contradiccion}
\left( 1+\frac{R^2}{4}  \right)^\alpha \frac{1}{R}\int_{R/2 < |x| <R} |u(x)|^2 dx \leq C.
\end{equation}

On the other hand, since $u$ is indeed a Herglotz wave function, from \eqref{limite_herglotz} we have that
\[\lim_{R\to\infty} \frac{1}{R}  \int_{R/2 < |x| <R} |u(x)|^2 dx=\frac{1}{2}\|u\|^2_{L},\]
which contradicts \eqref{contradiccion} when taking limits to infinity.
\end{remark}
\emph{Proof of Proposition \ref{prop9}.}  Any entire solution of the homogeneous Helmholtz equation $u$ admits the representation \eqref{expansionu}, and therefore, from 
\eqref{NuevaHW}, for any $R>0$ we have that
\[\frac{1}{R}\int_0^R\sum_{\ell=0}^\infty \sum_{j=1}^{\nu(\ell)}\left| a_{\ell j}\right|^2\left| J_{\mu (\ell)}(r) \right|^2 \frac{r\,dr}{(1+r^2)^\alpha}\leq C,\]
with $C$ an absolute constant. From here, using the dominated convergence theorem we get
\[
\sum_{\ell=0}^\infty \sum_{j=1}^{\nu(\ell)}
\left| a_{\ell j}\right|^2
\frac{1}{R}\int_0^R
\left| J_{\mu (\ell)}(r) \right|^2 \frac{r\,dr}{(1+r^2)^\alpha}
\leq C.
\]
And therefore, for $\ell=0,1,\ldots$, taking $R=2\mu(\ell)$, we have that
\[\sum_{j=1}^{\nu(\ell)}\left| a_{\ell j}\right|^2\frac{1}{\mu(\ell)}\int_{\mu(\ell)}^{2\mu(\ell)}\left| J_{\mu (\ell)}(r) \right|^2 \frac{r}{(1+r^2)^\alpha} dr \lesssim 1.\]
Taking into account that in the above integral $r\sim \mu(\ell)$ and that $\mu(\ell)=\ell+\frac{d-2}{2}\sim_d\ell$, using the estimate \eqref{bessel_uno} given below in the Appendix, we obtain that
\begin{equation}\label{pepino}
\sum_{j=1}^{\nu(\ell)}\frac{\left| a_{\ell j}\right|^2}{(1+\ell^2)^\alpha}\lesssim 1.
\end{equation}
Now, if we define 
\begin{equation}\label{definophi}
\phi(\theta)=\sum_{\ell=0}^\infty\sum_{j=1}^{\nu(\ell)}a_{\ell j }Y_\ell^j (\theta),
\end{equation}
from \eqref{def:H}, using \eqref{pepino}, we have that
$$
\|\phi\|_{H^{-\beta}(\mathbb{S}^{d-1})}^2=\sum_{\ell=0}^\infty\sum_{j=1}^{\nu(\ell)}\frac{|a_{\ell j}|^2}{(1+\ell^{\frac{1}{2}}(\ell+d-2)^{\frac{1}{2}})^{2\beta}}
\lesssim 
\sum_{\ell=0}^\infty\frac{1}{(1+\ell)^{2(\beta-\alpha)}}
<\infty,
$$
whenever $\beta>\alpha+1/2$.
And thus $\phi\in H^{-\beta}(\mathbb{S}^{d-1})$ for $\beta>\alpha +1/2$, and from \eqref{definophi} and \eqref{repW}, we have that $u=\mathcal{E}\phi$.
\hfill$\Box$

\emph{Proof of Corollary \ref{ultimo}.}
As we mentioned in Subsection \ref{preliminares}, the space of distributions in the sphere $C^\infty(\mathbb{S}^{d-1})^*$ can be identified with all the series given in \eqref{representationFunct} such that the coefficients $\widehat{f}_{\ell j}$ satisfy the growth condition given in \eqref{growth}.
From here, it is easy to see that 
\begin{equation*}
\distr=\underset{\alpha>0}{\bigcup} H^{-\alpha}(\mathbb{S}^{d-1}),
\end{equation*}
and thus, the result follows from Propositions \ref{donsalvador_1} and \ref{prop9}.\hfill$\Box$
\section{Appendix}
In this section we will prove Lemma \ref{lema_principal}, which involves the Bessel functions. The following results concerning these functions will be needed. They can be found in \cite{BRV} and \cite{BBR}.		
\begin{lemma}\cite[Lemma 1]{BRV} 
	\label{herramientas} 
	Let $\mu \geq 1/2$. 
	\begin{enumerate}
		\item 
		If $0<r<1$, then
		\begin{equation}\label{Bessel_cero}
		\left| J_\mu (r)\right| \lesssim \frac{1}{\Gamma (\mu +1)}\left( \frac{r}{2}  \right)^\mu.
		\end{equation}
		\item 
		If $1 \leq r \leq  \frac{\mu}{4}$ and we write $r=\mu\sech\alpha_\mu(r)$ with $ e^{\alpha_\mu (r)}=\frac{\mu}{r}+\frac{\sqrt{\mu^2-r^2}}{r} $, then
		$\tanh \alpha_\mu (r)=\frac{\sqrt{\mu^2-r^2}}{\mu} $ and
		\begin{equation}\label{Bessel_entre_1_y_n/2}
		\left| J_\mu ( \mu sech \alpha_\mu (r) ) \right| \lesssim 
		\frac{e^{\mu \tanh \alpha_\mu (r)}}{\mu^{1/2}e^{\mu \alpha_\mu(r)}} 
		\sim  \frac{e^{\sqrt{\mu^2-r^2}}}{\mu^{1/2}\left( \frac{\mu}{r}+\frac{\sqrt{\mu^2-r^2}}{r}   \right)^{\mu}}. 
		\end{equation}
		\item 
		If $r\ge 2\mu$, then
		\begin{equation}\label{granded}
		\left| J_\mu(r)   \right|  \lesssim r^{-1/2}.
		\end{equation}
	\end{enumerate}
\end{lemma}
\begin{lemma}\cite[Lemma 5]{BBR} \label{herramientas2}
	Let $\mu \geq 1/2$. If $a \geq 1 $, then
	\begin{equation}\label{bessel_uno}
	\int_{\mu/a}^{2 \mu} \left| J_\mu (r) \right|^2 dr \sim 1.
	\end{equation}
	Moreover, it holds
	\begin{equation}\label{bessel_uno_1}
1\lesssim(1+\mu)^2   \int_0^{\infty}     \left|   J_\mu (r) \right|^2 \frac{rdr}{<r>^3},
\end{equation}	
where $<r>=(1+r^2)^{1/2}.$	
\end{lemma}
\emph{Proof of Lemma \ref{lema_principal}}. 
It is enough to prove the following estimates:
\begin{equation}
\label{parte1}
\frac{1}{R}\int_{0}^R  \left|  J_\mu (r)  \right|^2\frac{r\,dr}{(1+r^2)^s}         
\lesssim  
(1+\mu)^{2-2s}
\int_0^\infty     \left|   J_\mu (r) \right|^2 \frac{r\,dr}{<r>^3},
\quad\forall R>0,
\end{equation}
\begin{equation}
\label{parte2}       
(1+\mu)^{2-2s}
\int_0^\infty     \left|   J_\mu (r) \right|^2 \frac{r\,dr}{<r>^3}
\lesssim
\sup_{R>0}
\frac{1}{R}\int_{0}^R  \left|  J_\mu (r)  \right|^2\frac{r\,dr}{(1+r^2)^s}.
\end{equation}
We start with the proof of \eqref{parte1}.
We will argue in a different way depending on the values of $R>0$, so we distinguish the following cases:

(1)
If $0<R<1$, we can use (\ref{Bessel_cero}) and \eqref{bessel_uno_1} to write
\begin{align}
\frac{1}{R}\int_0^R \left| J_\mu (r)  \right|^2\frac{r\,dr}{(1+r^2)^s}
&\lesssim \frac{1}{2^{2 \mu}\,\Gamma^2(\mu+1)}
\nonumber
\\
&\lesssim
(1+\mu)^{2-2s}
\int_0^\infty     \left|   J_\mu (r) \right|^2 \frac{r\,dr}{<r>^3}.
\label{launo} 
\end{align}
Notice that to get \eqref{launo} we have use that for any $s\ge 0$
$$
(1+\mu)^{2s-2}\lesssim 2^{2 \mu}\,\Gamma^2(\mu+1)
$$
whenever $\mu$ is large enough. This can be seen using
the Striling's formula for the gamma function (see \cite[pp.12]{LE}), that is,
\begin{equation}\label{stirling}
\Gamma (x)=\sqrt{2\pi} x^{x-\frac{1}{2}}e^{-x}\left(1+O\left(\frac{1}{x}\right)\right), \qquad  x>0.
\end{equation}

(2)
Let $1\le R<\mu/4$. From the previous case, it is enough to prove that
\begin{eqnarray}
\label{aprobar}
\frac{1}{R}\int_1^R \left| J_\mu (r)  \right|^2\frac{r\,dr}{(1+r^2)^s}
\lesssim 
(1+\mu)^{2-2s}
\int_0^\infty     \left|   J_\mu (r) \right|^2 \frac{r\,dr}{<r>^3}.
\end{eqnarray}
Using \eqref{Bessel_entre_1_y_n/2}, we have that
\begin{equation*}
\frac{1}{R}\int_1^R \left| J_\mu (r)  \right|^2\frac{r\,dr}{(1+r^2)^s}
\lesssim
\int_1^{\mu/4} 
\frac{e^{2\sqrt{\mu^2-r^2}}}{ \left(  \frac{\mu}{r}  \right)^{2 \mu} \left( 1+\frac{\sqrt{\mu^2-r^2}}{\mu}   \right)^{2 \mu}}
\,dr
\lesssim
\left(\frac{e}{4}\right)^{2\mu}.
\end{equation*}
Estimate \eqref{aprobar} follows from here using \eqref{bessel_uno_1}, since 
$$
\left(\frac{e}{4}\right)^{2\mu}\lesssim (1+\mu)^{2-2s}
$$
for any $s>0$ whenever $\mu$ is large enough.

(3) Let $\mu /4 \leq R < 2 \mu $. From the previous two cases, it is enough to prove 
\begin{equation*}
\label{aprobar2}
\frac{1}{R}\int_{\mu/4}^R \left| J_\mu (r)  \right|^2\frac{r\,dr}{(1+r^2)^s}
\lesssim 
(1+\mu)^{2-2s}
\int_0^\infty     \left|   J_\mu (r) \right|^2 \frac{r\,dr}{<r>^3}.
\end{equation*}
This estimate follows from \eqref{bessel_uno}  with $a=4$ and \eqref{bessel_uno_1} since 
\begin{equation}
\label{rmu}
(1+r^2)^{-s}\lesssim (1+\mu)^{-2s},\qquad s\ge0,\ r\ge\mu/4.
\end{equation}

(4) Let $2\mu \leq R$. From the previous three cases, it is enough to prove 
\begin{equation*}
\label{aprobar3}
\frac{1}{R}\int_{2\mu}^R \left| J_\mu (r)  \right|^2\frac{r\,dr}{(1+r^2)^s}
\lesssim 
(1+\mu)^{2-2s}
\int_0^\infty     \left|   J_\mu (r) \right|^2 \frac{r\,dr}{<r>^3}.
\end{equation*}
This estimate follows from \eqref{granded}, \eqref{rmu} and \eqref{bessel_uno_1}.

In order to prove \eqref{parte2}, we split the integral in the left-hand side into four integrals,
$$
\int_0^\infty=\int_0^1+\int_1^{\mu/4}+\int_{\mu/4}^{2\mu}+\int_{2\mu}^\infty,
$$
so that it is enough to prove the result for each one of them.

For the first integral, using \eqref{Bessel_cero}, the Stirling's formula given in \eqref{stirling} and \eqref{bessel_uno} with $a=1$, we can write
\begin{align*}
(1+\mu)^{2-2s}
\int_0^1
\left|   J_\mu (r) \right|^2 \frac{r\,dr}{<r>^3}
&\lesssim
\frac{(1+\mu)^{2-2s}}{2^{2\mu}\,\Gamma^2(\mu+1)}
\\
&\lesssim
\frac{1}{(1+\mu)^{2s}}
\int_\mu^{2\mu}
\left|   J_\mu (r) \right|^2 \,dr
\\
&\sim
\frac{1}{2\mu}
\int_\mu^{2\mu}
\left|   J_\mu (r) \right|^2 \frac{r\,dr}{(1+r^2)^s}
\\
&\le
\sup_{R>0}
\frac{1}{R}\int_{0}^R  \left|  J_\mu (r)  \right|^2\frac{r\,dr}{(1+r^2)^s}.
\end{align*}

For the second integral, the argument is similar. Using \eqref{Bessel_entre_1_y_n/2} and \eqref{bessel_uno} with $a=1$, arguing as before, we get
\begin{align*}
(1+\mu)^{2-2s}
\int_1^{\mu/4}
\left|   J_\mu (r) \right|^2 \frac{r\,dr}{<r>^3}
&\lesssim
(1+\mu)^{2-2s}\left(\frac{e}{4}\right)^{2\mu}\mu
\\
&\lesssim
\frac{1}{(1+\mu)^{2s}}
\int_\mu^{2\mu}
\left|   J_\mu (r) \right|^2 \,dr
\\
&\le
\sup_{R>0}
\frac{1}{R}\int_{0}^R  \left|  J_\mu (r)  \right|^2\frac{r\,dr}{(1+r^2)^s}.
\end{align*}

In the third integral, since we are considering $r\sim \mu$, we have that
\begin{align*}
(1+\mu)^{2-2s}
\int_{\mu/4}^{2\mu}
\left|   J_\mu (r) \right|^2 \frac{r\,dr}{<r>^3}
&\lesssim
\frac{(1+\mu)^2}{\mu^3}
\int_{\mu/4}^{2\mu}
\left|   J_\mu (r) \right|^2 
\frac{r\,dr}{(1+r^2)^s}
\\
&\le
\sup_{R>0}
\frac{1}{R}\int_{0}^R  \left|  J_\mu (r)  \right|^2
\frac{r\,dr}{(1+r^2)^s}.
\end{align*}

Finally, for the fourth integral we use (\ref{granded}) and (\ref{bessel_uno}) with $a=1$ to write
\begin{align*}
(1+\mu)^{2-2s}
\int_{2\mu}^\infty
\left|   J_\mu (r) \right|^2 \frac{r\,dr}{<r>^3}
&\lesssim
\frac{(1+\mu)^{2-2s}}{\mu^2}
\\
&\lesssim
\frac{1}{(1+\mu)^{2s}}
\int_\mu^{2\mu}
\left|   J_\mu (r) \right|^2 \,dr
\\
&\le
\sup_{R>0}
\frac{1}{R}\int_{0}^R  \left|  J_\mu (r)  \right|^2\frac{r\,dr}{(1+r^2)^s}.
\end{align*}
\begin{flushright}
	$\Box$
\end{flushright}


J. A. Barcel\'o,  Departamento de Matem\'atica e Inform\'atica aplicadas a las Ingenier\'ias Civil y Naval, Universidad Polit\'ecnica de Madrid,  28040 Madrid, Spain. \\email: juanantonio.barcelo@upm.es.

M. Folch-Gabayet, Instituto de Matem\'{a}ticas,
Universidad Nacional Aut\'{o}noma de M\'{e}xico, Ciudad Universitaria,
Ciudad de M\'{e}xico, 04510, M\'{e}xico. \\email: folchgab@matem.unam.mx.

T. Luque, Departamento de An\'alisis Matem\'atico y Matem\'atica Aplicada, Universidad Complutense de Madrid, 28040 Madrid, Spain. \\email: t.luque@ucm.es.

S. P\'erez-Esteva, Instituto de Matem\'aticas, Unidad de Cuernavaca, Universidad Nacional Aut\'onoma de M\'exico, M\'exico. \\email: spesteva@im.unam.mx.

M. C. Vilela, ETSI Navales, Departamento de Matem\'atica e Inform\'atica aplicadas a las Ingenier\'ias Civil y Naval,
Universidad Polit\'ecnica de Madrid, 28040 Madrid, Spain.\\
email:  maricruz.vilela@upm.es.

	\end{document}